\documentclass[11pt,a4paper]{article}
\usepackage{graphicx} % Required for inserting images
\usepackage{amsfonts,amsmath,amsthm,amssymb}
\usepackage{xcolor}
\usepackage{geometry}[margin=1in]
\newtheorem{theorem}{Theorem}[section]

\newtheorem{remark}[theorem]{Remark}
\newtheorem{example}[theorem]{Example}
\newtheorem{hypothesis}{Hypothesis}

\usepackage{hyperref}
\providecommand{\keywords}[1]
{
  \small	
  \textbf{\textit{Keywords---}} #1
}
\newcommand{\footremember}[2]{%
    \footnote{#2}
    \newcounter{#1}
    \setcounter{#1}{\value{footnote}}%
}

\title{Stochastic Solutions for Hyperbolic PDE}

\author{%
  Abdol-Reza Mansouri\footremember{alley}{Queen's University. mansouri@queensu.ca}%
  \and Zachary Selk\footremember{trailer}{Queen's University. zachary.selk@queensu.ca}%
  }

\date{\today}

\begin{document}

\maketitle

\begin{abstract}
    The theory of stochastic representations of solutions to elliptic and parabolic PDE has been extensive. However, the theory for hyperbolic PDE is notably lacking. In this short note we give a stochastic representation for solutions of hyperbolic PDE. 
\end{abstract}


\section{Introduction}
There has been extensive work done on representing solutions of parabolic and elliptic PDE stochastically. For example, one simplified version of the celebrated Feynman-Kac theorem states that the heat equation
$$\begin{cases}
     \frac{\partial u }{\partial t}(t,x)=\frac{1}{2} \Delta u(t,x)\text{ on }(0,\infty)\times \mathbb R^d,\\
    u(0,x)=f(x)\text{ on }\{0\}\times \mathbb R^d,
\end{cases}$$
with $f$ continuous and imposing sufficient decay, has solution $E[f(x+B(t))]$, where $B$ is a standard $d$-dimensional Brownian motion. 

Additionally, for domain $D\subset \mathbb R^d$ Kakutani (see \cite{MR0014646,MR0014647}) showed that the Dirichlet problem
$$\begin{cases}
    \Delta u =0 \text{ in } D\\
    u\big|_{\partial D}=f
\end{cases}$$
where $f$ is a continuous function defined on the boundary $\partial D$ has solution $E[f(x+B(\tau))]$ where $\tau=\inf\{t:x+B(t)\in \partial D\}$ is the exit time of Brownian motion from the domain $D$. 

However, for hyperbolic PDE such as the wave equation, the literature has been lacking. In the preprint \cite{Chaterjee}, the author presents a way of generating solutions to hyperbolic PDE by making joint use of the exit time analysis for elliptic PDE and the specific distribution of Cauchy random variables, thereby converting solutions of elliptic PDE to solutions of hyperbolic PDE. In particular, the author shows the following theorem. 

\begin{theorem}
    Let $d\geq 1$ and let $D$ be a bounded open connected subset of $\mathbb R^d$ with boundary $\partial D$. Let $f:\mathbb R\times \partial D\to \mathbb R$ be any bounded measurable function. Let $B$ be a $d$-dimensional standard Brownian motion, let $Z$ be a standard Normal random variable and let $X$ be a standard Cauchy random variable. Denote by  $\tau$ the exit time of $B$ from set $D$. Then $$u(t,x)=E[f(tX+\sqrt{\tau}Z, x+ B(\tau))]$$ is twice continuously differentiable and satisfies the wave equation $\frac{\partial^2 u}{\partial t^2}=\Delta u$.
\end{theorem}

The idea behind this representation is that the function $v$ given by $v(t,x):=E[f(t+\sqrt{\tau}Z,B(t)+x)]$ solves the $d+1$ dimensional space-time Laplace equation. The Cauchy random variable then changes the sign of the second time derivative. Although innovative and generating a broad class of solutions stochastically, the author comments that this representation has several drawbacks. In particular the limitations are:

\begin{enumerate}
    \item This only generates bounded solutions.
    \item Even among bounded solutions, this only generates solutions for which $\lim_{t\to\infty}\partial_t u(t,x)=0$.
    \item The boundary data is ``time evolution on the boundary" which might be seen as unnatural, and further the relationship with this boundary data and more common forms is unclear.
    \item The Cauchy random variable introduces some randomness that might be seen as unnatural.
    \item This method provides a representation for solutions only to the wave equation.
\end{enumerate}

We also refer to \cite{Mueller-Wave,Goldstein-Hyperbolic,Wang-wave} for related approaches. In this note, we use complex analysis and a version of ``Wick rotation" (see \cite{Wick-rotation}) to give a new representation. We assume the following hypotheses. 

\begin{hypothesis}\label{hyp:L}Let $d\geq 1$ where $d\in \mathbb N$. Let $b:\mathbb R^d\to \mathbb R^d, \sigma:\mathbb R^d\to \mathbb R^{d\times d}$ be Lipschitz, thus implying that the stochastic differential equation 
    $$dX^{x_0}(t)=b(X^{x_0}(t))dt+\sigma(X^{x_0}(t))dB(t)$$
    with $X^{x_0}(0)=x_0$ has a unique strong solution, and its infinitesimal generator is given by
    $$L\big|_x=\sum_{i=1}^d b_i(x)\frac{\partial }{\partial x_{i}}+\sum_{i,j=1}^d a_{i,j}(x)\frac{\partial^2}{\partial x_i \partial x_j},$$
    where $a(x)=\frac{1}{2}\sigma(x)\sigma^T(x)$, with $x\in \mathbb R^d$.
\end{hypothesis}
\begin{hypothesis}\label{hyp:domain}
    Let $D\subset \mathbb R^d$ be a bounded open connected set with regular boundary - that is, the exit time $\tau_D:=\inf\{s>0:X^{x_0}(s)\in \partial D\}$ is $0$ with probability $1$ if $x\in \partial D$. 
\end{hypothesis}
One sufficient condition for Hypothesis \ref{hyp:domain} to hold is that the boundary $\partial D$ is $C^1$. In what follows we denote by $C^2( \mathbb C\times \bar D, \mathbb C)$ the set of continuous functions $h: \mathbb C\times \bar D\to \mathbb C$ that are complex differentiable on $\mathbb C$ in the first component and twice real differentiable on $D$ in the second component. 
\begin{hypothesis}\label{hyp:f}
Let $Z$ be a standard normal random variable independent of $X$ and let $B$ be a standard Brownian motion independent of $X$. We define the space-time domain $D'=(0,\infty)\times D$ with boundary $\partial D'=([0,\infty) \times \partial D)\cup (\{0\}\times D)$. For each $(z,x_0)\in (\mathbb C\setminus \{0\})\times \bar D$ we define $\tau=\inf\{s>0:(z+B(s),X^{x_0}(s))\in \partial ((\mathbb C\setminus\{0\})\times D)\}$, i.e. $\tau$ is the first exit time of the process $(z+B,X^{x_0})$ from the set $(\mathbb C \setminus\{0\})\times D$. Let $f\in C^2( \mathbb C\times \bar D, \mathbb C)$ be so that $$E[|f(z+i\sqrt \tau Z, X^{x_0}(\tau))|]<\infty$$
for all $(z,x_0)\in \mathbb C\times \bar D.$ Suppose for all compact sets $G\subset \mathbb C$ there is a $g_G:\mathbb C\times \bar D\to [0,\infty)$, with $E[g_G(i \sqrt \tau Z, X^{x_0}(\tau))]<\infty,$ such that
$$|\partial_t f(z+i\sqrt \tau Z, X^{x_0}(\tau))|\leq g_G(i \tau Z, X^{x_0}(\tau))$$
for all $z\in G$.
\end{hypothesis}
\begin{remark}
    Note that $\tau$ is defined as the exit time of the process $(z+B,X^{x_0})$ from the ``partially complexified" space-time domain  $ (\mathbb C\setminus\{0\})\times D$ in anticipation of the complexification of the first independent variable. Note also that the stopping time $\tau$ is defined as the exit time from $(\mathbb C\setminus \{0\})\times D$ and not from $\mathbb C \times D$; this allows us to recover the initial conditions in addition to the boundary conditions for the hyperbolic problem. 
\end{remark}
Our theorem is the following. 

\begin{theorem}\label{Theorem:Main}
Assume Hypotheses \ref{hyp:L}, \ref{hyp:domain}, \ref{hyp:f} hold. Then one solution to the equation
    \begin{equation}
        \begin{cases}
            \frac{1}{2}\frac{\partial^2 u}{\partial t^2}&=L u,\\
            u\big|_{\partial D'}&=f,\\
            \partial_t u\big|_{\{0\}\times  D}&=\partial_t f(0,x),
        \end{cases}
    \end{equation}
    is given by
    \begin{equation}
    u(t,x)=E[  f(t+i\sqrt{\tau} Z,X^x(\tau))]
    \end{equation}
    for $(t,x)\in [0,\infty)\times \bar D$ where $\tau$ is the stopping time associated to $(z,x_0)=(it,x)\in \mathbb C\times \bar D.$
\end{theorem}

The idea behind this representation is to apply a ``Wick rotation". Wick rotation consists in ``rotating" time by $\pi/2$ by sending $t\mapsto it$. If $v(t,x)$ solves a $d+1$ dimensional space-time elliptic PDE then $(t,x)\mapsto u(t,x):=v(it,x)$ solves a $d$ space and $1$ time dimensional hyperbolic PDE. The approach in \cite{Chaterjee} uses a Cauchy random variable to accomplish this, but using the imaginary $i$ instead allows us to relax the requirements on $f$.

Our result overcomes some of the limitations of the approach in proposed in \cite{Chaterjee}. We have no boundedness constraints on $f$, nor do we require that $\lim_{t\to\infty} \partial_t u(t,x)=0$. Furthermore, our representation allows for more natural initial and boundary conditions. Our representation additionally allows for more general hyperbolic PDE than the classical wave equation. The only limitation of our approach is that it does not cover all initial and boundary conditions - just ones that can be ``unified" in a single function $f$.

\begin{remark}
    Given any $f:\mathbb C\times \bar D\to \mathbb C$  and $\phi: \bar D\to \mathbb R$ with $\phi(x)=0$ if $x\in \partial D$ we can consider $ f_1(z,x):=f(z,x)+z\phi(x)$. In this case, the initial velocity is $\partial_t f(0,x)+\phi(x)$ and the initial position is still $f(0,x).$ Therefore the the initial velocity is not actually constrained by the initial position, contrary to what the statement of Theorem \ref{Theorem:Main} might suggest.
\end{remark}

One strength of stochastic representations of solutions to PDE is they readily yield numerical Monte-Carlo methods. However, as in the example given in \cite{Chaterjee}, it is still possible to obtain explicit solutions in some cases.

\begin{example}
Let $D=(0,1)\subset \mathbb R$ and let $L=\frac{1}{2}\frac{\partial^2}{\partial x^2}$. Consider the hyperbolic PDE
\begin{equation}
        \begin{cases}
           \frac{1}{2}\frac{\partial^2 u}{\partial t^2}&=\frac{1}{2}\frac{\partial^2 u}{\partial x^2}\\
            u(t,0)&=e^{t}\\
            u(t,1)&=e^{t+1}\\
            \partial_t u(0,x)&=e^{x}
        \end{cases}
\end{equation}
Letting $f(z,x)=e^{x+z}$ agrees with the boundary and initial conditions. Then
    \begin{align*}
        E[  f(t+i\sqrt{\tau} Z,X^x(\tau))]&=E[  f(t+i\sqrt{\tau} Z,x+B(\tau))]\\
        &=E[E[  f(t+i\sqrt{\tau} Z,x+B(\tau))|\sigma(B(\tau))]]\\
        &=E[E[e^{x+B(\tau)}e^{t+i\sqrt{\tau} Z}|\sigma(B(\tau))]]\\
        &=E[e^{x+B(\tau)}E[e^{t+i\sqrt{\tau} Z}|\sigma(B(\tau))]]\\
        &=e^{x+t} E[e^{B(\tau)-\tau/2}].
    \end{align*}
As $(e^{B(t)-t/2})_{t\geq 0}$ is a martingale and $\left|e^{B(t\wedge \tau)-(t\wedge \tau)/2}\right|\leq e^1$ a.s. for all $t\geq 0$, Doob's Optional Stopping Theorem yields $E[e^{B(\tau)-\tau/2}]=1$. Therefore the solution $u$ is given by $u(t,x)=e^{x+t}$ for all $(t,x)\in [0,\infty)\times [0,1]$. We emphasize that this solution is unbounded, in contrast to the solutions obtainable in \cite{Chaterjee}.
\end{example}

\section{Proof of Theorem \ref{Theorem:Main}}
\begin{proof}
The proposed solution $u$ satisfies
\begin{align*}
    u(t,x)&=E[f(t+i\sqrt{\tau} Z,X^x(\tau))]\\
    &=E[  f((it-\sqrt{\tau} Z)/i,X^x(\tau))]\\
    &:=E[\hat f(it-\sqrt{\tau}Z,X^x(\tau))],
\end{align*}
where $\hat f(z,x):=f(z/i,x)$. Note that $\hat f\in C^2(\mathbb C \times \bar D, \mathbb C)$ as well. Additionally, $\hat f$ enjoys the same integrability conditions as outlined in Hypothesis \ref{hyp:f}. Furthermore, as $Z$ is standard normal and independent of $X$, $u$ can be written as
\begin{align*}
    u(t,x)&=E[\hat f(it+B(\tau),X^x(\tau))],
\end{align*}
where $B$ is a standard Brownian motion independent of $X.$

Define for $z\in (\mathbb C\setminus\{0\})$ and $x\in  D$ the function $v(z,x)=E[ \hat f(z+B(\tau), X^x(\tau))]$ where we recall that $\tau=\inf\{s>0:(z+B(s),X^x(s))\in \partial ((\mathbb C\setminus\{0\})\times D)\}$. We will denote this by
$$v(z,x):=E^{(z,x)}[ \hat f(\Tilde{X}(\tau))],$$
where $\Tilde{X}$ is the $d+1$ dimensional space-time process $(B,X)$. The notation $E^{(z,x)}$ denotes the expectation of $\tilde X$ started at the space-time point $(z,x)\in \mathbb C\times D$. As $X$ is a Feller process, so is $\tilde X$ with generator $\frac{1}{2}\frac{\partial^2}{\partial z^2}+L$, where $L$ acts on the $x$ variable only. As $\hat f\in C^2(\mathbb C\times \bar D, \mathbb C)$, Lemma 9.2.3 in \cite{Oksendal} along with holomorphicity (recall that if $h:\mathbb C\to \mathbb C$ is holomorphic, with $h(\xi+i\eta)=h_1(\xi,\eta)+ih_2(\xi,\eta)$ where $\xi,\eta\in \mathbb R$, and $z=\xi+i\eta$, then $\frac{d}{dz}h(z)=\frac{\partial}{\partial \xi}h_1(\xi,\eta)+i\frac{\partial}{\partial \xi}h_2(\xi,\eta)$) yield that $v$ solves the equation
\begin{equation}
    \frac{1}{2}\frac{\partial^2 v}{\partial z^2}+Lv=0.
\end{equation}
As $u(t,x)=v(it,x)$, the chain rule and holomorphicity in the first variable then imply 
\begin{equation}
   - \frac{1}{2}\frac{\partial^2 u}{\partial t^2}+Lu=0.
\end{equation}
The only thing to check now are the boundary conditions. First, if $(t,x)\in \partial D'$ then Hypothesis \ref{hyp:domain} implies that $\tau=0$ and $u(t,x)=f(t,x)$. 

Second, by Hypothesis \ref{hyp:f} we may apply the dominated convergence theorem to interchange differentiation and expectation to get
\begin{equation}
    \partial_t u(t,x)=E[\partial_t f(t+i\sqrt{\tau} Z,X^x(\tau))].
\end{equation}
If $t=0$ and $x\in D$ then $\tau=0$ and 
\begin{equation}
    \partial_t u(0,x)=E[\partial_t f(0,x)]=\partial_t f(0,x).
\end{equation}
\end{proof}
\section*{Acknowledgments}
The authors gratefully acknowledge the anonymous referee whose comments have greatly improved the presentation of this paper. This work was supported in part by the Natural Sciences and Engineering Research Council of Canada.

\bibliographystyle{plain}
\bibliography{bibliography}
\end{document}